\documentclass[11pt]{article}
\usepackage{amssymb}
\usepackage{amsmath}
\usepackage{amsthm}

 \usepackage{color}
 \usepackage{ulem}

\numberwithin{equation}{section}

\date{}

\title{\bf  Notes on solution maps of abstract FDEs}

\author{Xiao-Qiang Zhao\thanks{Research supported in part by the NSERC of Canada.}\\
Department of Mathematics and Statistics\\
Memorial University of Newfoundland\\
St. John's, NL A1C  5S7, Canada\\
E-mail:\,  zhao@mun.ca}

\begin{document}
\maketitle


Let $\tau$ be a positive real number, $X$ be a Banach space,
and  $C:=C([-\tau,0],X)$.  For any  $\phi\in C$, define
$\|\phi\|=\sup\limits_{-\tau \leq
	\theta \leq 0}\|\phi(\theta)\|_X$.
Then $(C, \|\cdot\|)$ is a Banach space.
Let $A$ be the infinitesimal generator of a $C_0$-semigroup 
$\{T(t)\}_{t\geq 0}$ on $X$.  Assume that $T(t)$ is compact
for each $t>0$, and there exists $M>0$ such that $\|T(t)\|\leq M$ for all $t\geq 0$. 

We consider the  abstract functional differential equation
\begin{equation}
\label{FDE}
\begin{split}
& \frac{du(t)}{dt}=Au(t)+F(t,u_t), ~~~~t>0, \\
& u_0=\phi\in C.
\end{split}
\end{equation}
Here $F: [0,\infty)\times C\rightarrow X$ is continuous and maps bounded sets into bounded sets,  and $u_t\in C$ is defined by
$u_t(\theta)=u(t+\theta), \,  \forall \theta\in [-\tau,0]$.

\

\noindent
{\bf Theorem  A. } Assume that for each $\phi\in C$,
equation \eqref{FDE} has a unique solution $u(t,\phi)$ on
$[0,\infty)$, and  solutions of \eqref{FDE} are
uniformly bounded in the sense that for any bounded subset $B_0$
of $C$, there exists a bounded subset $B_1=B_1(B_0)$ of $C$ such that
$u_t(\phi)\in B_1$ for all $\phi\in B_0$ and $t\geq 0$.
Then for any given  $r>0$, there exists an
equivalent norm $\|\cdot\|_r^*$ on $C$ such that 
the solution maps
$Q(t):=u_t$ of equation \eqref{FDE} satisfy
$\alpha(Q(t)B)\leq e^{-rt}\alpha(B)$
for any bounded subset $B$ of  $C$ and $t\geq 0$, where $\alpha$ is the Kuratowski  measure of noncompactness in $(C,\|\cdot\|_r^*)$.

\

\noindent{\it Proof.}  Define $\|x\|^*=\sup_{t\geq 0}\|T(t)x\|, \,  \, \forall x\in X$. Then $\|x\|\leq \|x\|^*\leq M \|x\|$,
and hence,  $\|x\|^*$ is an equivalent norm on $X$. It is easy to see that
\[
\|T(t)x\|^*=\sup_{s\geq 0}\|T(s)T(t)x\|=
\sup_{s\geq 0}\|T(s+t)x\|\leq \|x\|^*, \, \, \,  \forall x\in X,\, \, t\geq 0,
\]
which implies that $\|T(t)\|^*\leq 1$ for all $t\geq 0$.
Thus, without loss of generality, we assume that  $M=1$. 

Let $r>0$ be given.  Note that for each $\phi\in C$, the solution $u(t,\phi)$ of \eqref{FDE} satisfies the following integral equation
\begin{equation}
\label{xq1}
\begin{split}
& u(t)=\hat T(t)\phi(0)+\int_0^t \hat T(t-s)\hat F(s,u_s)ds, ~~~~t\geq 0, \\
& u_0=\phi\in C,
\end{split}
\end{equation}
where  $\hat T(t)=e^{-rt}T(t)$ and $\hat F(t,\varphi)=r\varphi(0)+F(t,\varphi), \, \, \forall t\geq 0, \, \varphi\in C$.
Then $\hat T(t)$ is also a $C_0$-semigroup on $X$ and
$\|\hat T(t)\|\leq e^{-rt},\, \, \forall t\geq 0$.
Let $h(\theta)=e^{-r\theta}, \forall~\theta\in
[-\tau,0]$, and define
\[\|\phi\|_r^*=\sup\limits_{-\tau\leq \theta\leq 0}\frac{\|\phi(\theta)\|_X}{h(\theta)},~~\forall~\phi\in C.\]
Then $\frac{1}{h(-\tau)}\|\phi\|_C\leq \|\phi\|_r^*\leq
\|\phi\|_C$, and hence $\|\cdot\|_r^*$ is equivalent to
$\|\cdot\|_C$.  Clearly,  $\|\phi(0)\|_X\leq  \|\phi\|_r^*, \, \forall \phi\in C$.  Define
\begin{align*}
(L(t)\phi)(\theta)=\left\{
\begin{array}{ll}
\hat T(t+\theta)\phi(0), ~~& t+\theta>0,\\
\phi(t+\theta), ~~& t+\theta\leq 0,
\end{array}
\right.
\end{align*}
and
\begin{align*}
(\bar{Q}(t)\phi)(\theta)=\left\{
\begin{array}{ll}
\int_0^{t+\theta} \hat T(t+\theta-s)\hat F(s,u_s(\phi))ds, & t+\theta>0,\\
0, ~~& t+\theta\leq 0.
\end{array}
\right.
\end{align*}
Thus, $Q(t)\phi=L(t)\phi+\bar{Q}(t)\phi, \forall~t\geq
0, \phi\in C,$ that is, $Q(t)=L(t)+\bar{Q}(t), \forall~t\geq 0$.

We first show that $L(t)$ is an $\alpha$-contraction on $(C,\|\cdot\|_r^*)$ for each  $t>0$. It is easy to see that 
$L(t)$ is compact for each $t>\tau$.  Without loss of 
generality, we may assume that  $t\in (0,\tau]$ is fixed.
For any $\phi\in C$, we have 
\begin{align*}
\|L(t)\phi\|_r^* &=\sup\limits_{-\tau\leq \theta\leq 0}
\frac{\|(L(t)\phi)\|_X}{h(\theta)}\\
& \leq \max\left\{ \sup\limits_{-\tau\leq \theta\leq
	-t}\frac{\|\phi(t+\theta)\|_X}{h(t+\theta)}\frac{h(t+\theta)}{h(\theta)},
~~\sup\limits_{-t\leq \theta\leq
	0}\frac{\|\hat T(t+\theta)\phi(0)\|_X}{h(\theta)} \right\}\\
& \leq \max\left\{ e^{-rt}\|\phi\|_r^*, ~~\sup\limits_{-t\leq
	\theta\leq
	0}\frac{ e^{-r(t+\theta)}\|\phi(0)\|_X }{ h(\theta) } \right\}\\
& = \max\left\{ e^{-rt}\|\phi\|_r^*, \, \,  e^{-rt}\|\phi(0)\|_X\right\} \leq e^{-rt}\|\phi\|_r^*,
\end{align*}
which implies that 
$\alpha(L(t)B)\leq e^{-rt}\alpha(B)$
for any bounded subset $B$ of $C$. Thus, this contraction
property holds true for all $t> 0$.

Next we prove that $\bar{Q}(t): C\rightarrow C$ is compact for
each $t>0$.  Let  $t>0$ and the bounded subset $B$ 
of $C$ be given.  By the uniform boundedness of solutions,  there exists a real number $K>0$ such that 
$\|\hat F(s,u_s(\phi))\|_X\leq K,\, \, \forall s\in [0,t], \, \phi\in B$.
It then follows that $\bar{Q}(t)B$ is bounded in $C$. 
We only need to show that  $\bar{Q}(t)B$ is precompact in $C$.   In view of the 
Arzela--Ascoli theorem for the space $C:=C([-\tau,0],X)$, it suffices to
prove that  (i) for each $\theta\in [-\tau,0]$,  the set $\{(\bar Q(t)\phi)(\theta): \,  \phi\in B\}$ is precompact in $X$; and  (ii) the set $\bar{Q}(t)B$ is  equi-continuous in $\theta\in [-\tau, 0]$.  
	Clearly,  statement (i)  holds true if $t+\theta\leq 0$. In the case where
	$t+\theta>0$, for any given $\epsilon \in (0,t+\theta)$, we have
	\begin{eqnarray*}
		(\bar{Q}(t)\phi)(\theta)&=&
		\int_0^{t+\theta-\epsilon} \hat T(t+\theta-s)\hat F(s,u_s(\phi))ds
		+\int_{t+\theta-\epsilon}^{t+\theta} \hat T(t+\theta-s)\hat F(s,u_s(\phi))ds\\
		&=&	\hat T(\epsilon)\int_0^{t+\theta-\epsilon} \hat T(t+\theta-\epsilon-s)\hat F(s,u_s(\phi))ds\\
		&\,&\, 
		+	\int_{t+\theta-\epsilon}^{t+\theta} \hat T(t+\theta-s)\hat F(s,u_s(\phi))ds.
	\end{eqnarray*}
	Define
	\[
	S_1:= \left\{\hat T(\epsilon)\int_0^{t+\theta-\epsilon} \hat T(t+\theta-\epsilon-s)\hat F(s,u_s(\phi))ds: \,  \phi\in B\right\}
	\]
	and
	\[
	S_2:=\left\{\int_{t+\theta-\epsilon}^{t+\theta} \hat T(t+\theta-s)\hat F(s,u_s(\phi))ds: \,  \phi\in B\right\}.
	\]
	Let $\hat \alpha$ be the Kuratowski measure of noncompactness in $X$.
	Since $\hat T(\epsilon)$ is compact, it follows that
	\[
	\hat \alpha\left(\{(\bar Q(t)\phi)(\theta): \,  \phi\in B\}\right)
	\leq \hat \alpha(S_1)+\hat \alpha(S_2)\leq 0+2K\epsilon=2K\epsilon.
	\]
	Letting $\epsilon\to 0^+$, we obtain 
	$\hat \alpha \left(\{(\bar Q(t)\phi)(\theta): \,  \phi\in B\}\right)=0$, which implies that 
	the set $\{(\bar Q(t)\phi)(\theta): \,  \phi\in B\}$ is precompact in $X$.
It remains to verify  statement (ii).
Since $\hat T(s)$ is compact for each $s>0$, $\hat T(s)$ is continuous in the
uniform operator topology for $s>0$ (see \cite[Theorem 2.3.2]{Pazy}).  It then follows that for any $\epsilon\in (0,t)$, there 
exists a $\delta=\delta (\epsilon)<\epsilon$ such that
\begin{equation}\label{norm-continuity}
\|\hat T(s_1)-\hat T(s_2)\|< \epsilon, \quad  \forall  s_1, s_2\in [\epsilon,t]
\, \,  \text{with}\, \,  |s_1-s_2|<\delta.
\end{equation}
We first consider the case where $t\in (0,\tau]$.  It is easy to see that 
\begin{equation}\label{equiE1}
\|(\bar{Q}(t)\phi)(\theta)\|_X\leq K(t+\theta)\leq K\epsilon,
\, \,  \forall \theta\in [-t,-t+\epsilon], \, \phi\in B.
\end{equation}
For any $\phi\in B$  and $\theta_1, \theta_2\in [-t+\epsilon,0]$
with $0<\theta_2-\theta_1<\delta$, it follows from \eqref{norm-continuity} that
\begin{eqnarray}\label{equiE2}
&\,&\left\|(\bar{Q}(t)\phi)(\theta_2)-(\bar{Q}(t)\phi)(\theta_1)\right\|_X 
\nonumber\\
&\,&=\left\|\int_0^{t-\epsilon+\theta_1}(\hat T(t+\theta_2-s)-\hat T(t+\theta_1-s))
\hat F(s,u_s(\phi))ds\right\|_X   \nonumber\\
&\, &\quad +\left\|\int_{t-\epsilon+\theta_1}^{t+\theta_2}
\hat T(t+\theta_2-s)\hat F(s,u_s(\phi))ds\right\|_X   \nonumber\\
&\,&\quad +\left\|-\int_{t-\epsilon+\theta_1}^{t+\theta_1}
\hat T(t+\theta_1-s)\hat F(s,u_s(\phi))ds\right\|_X   \nonumber\\
&\,&\leq \epsilon Kt+K(\theta_2-\theta_1+\epsilon)+K\epsilon  \nonumber\\
&\,&< (t+3)K\epsilon.
\end{eqnarray}
Combining \eqref{equiE1} and \eqref{equiE2},  we then obtain
\[
\|(\bar{Q}(t)\phi)(\theta_2)-(\bar{Q}(t)\phi)(\theta_1)\|_X<
2K\epsilon+(t+3)K\epsilon=(t+5)K\epsilon,
\]
for all $\theta_1, \theta_2\in [-t,0]$  with  $0\leq \theta_2-\theta_1<\delta$.
Since  $(\bar{Q}(t)\phi)(\theta)=0,  \forall \theta\in [-\tau, -t]$, it
follows  that 
$\bar{Q}(t)B$ is  equi-continuous in $\theta\in [-\tau, 0]$.
In the case where $t>\tau$, for any $\epsilon\in (0,t-\tau)$,
the  estimate in  \eqref{equiE2} implies that 
\[
\|(\bar{Q}(t)\phi)(\theta_2)-(\bar{Q}(t)\phi)(\theta_1)\|_X<
(t+3)K\epsilon,
\]
for all $\theta_1, \theta_2\in [-\tau,0]$  with  $0\leq \theta_2-\theta_1<\delta$,  and $\phi\in B$. Thus, 
$\bar{Q}(t)B$ is  equi-continuous in $\theta\in [-\tau, 0]$.
It then follows that  $\bar{Q}(t): C\rightarrow C$ is compact for
each $t>0$.

Consequently, for any $t>0$  and any bounded subset $B$ of $C$, we have
\[
\alpha(Q(t)B) \leq \alpha(L(t)B)+\alpha(\bar{Q}(t)B)\\
\leq e^{-rt}\alpha(B).
\]
This completes the proof.
\

As an application example, we consider the following $\omega$-periodic reaction-diffusion system
\begin{align}\label{RDE}
\left\{\begin{array}{ll}
\frac{\partial u}{\partial t}=D\triangle
u+f(t, u_t), ~~& x\in\Omega,\, \,  t>0,\\
\frac{\partial u}{\partial \nu}=0,~~& x\in\partial \Omega, \, \, t>0,
\end{array}
\right.  
\end{align}
where $D=diag\left(d_1,\dots, d_m\right)$ with each $d_i>0$, 
$\Omega\subset \mathbb{R}^n$ is a bounded domain with the smooth boundary $\partial \Omega$, 
and $f(t,\phi)$ is $\omega$-periodic in $t\in [0,\infty)$ for some
$\omega>0$. 
Let $Y:=C(\bar{\Omega},\mathbb{R}^m)$ and assume that 
$f$ is continuous and maps bounded subsets of $[0,\infty)\times C([-\tau,0],Y)$ 
into bounded subsets of $Y$.
Let  $T(t)$ be the semigroup on $Y$ generated by
$\frac{\partial u(t,x) }{\partial t}=D\triangle u(t,x)$
subject to the boundary condition $\frac{\partial u}{\partial \nu}=0$.
It is easy to see that $\|T(t)\|\leq 1, \forall t\geq 0$.
By Theorem A, we then have the following result.

\

\noindent
{\bf Theorem  B. } Assume that solutions of system \eqref{RDE}
exist uniquely on $[0,\infty)$ for any initial data in  $C:=C([-\tau,0],Y)$ and are uniformly bounded. Then for each $r>0$, there exists an
equivalent norm $\|\cdot\|_r^*$ on $C$ such that for each $t>0$,
the solution map $Q(t)=u_t$ of system \eqref{RDE} is 
an $\alpha$-contraction on $(C,\|\cdot\|_r^*)$ with 
the contraction constant being $e^{-r t}$. 

\

\noindent 
{\bf Remark.}  By using the theory of evolution operators 
(see, e.g., \cite[Section II.11]{Hessbook} and \cite[Section 5.6]{Pazy}),
one may extend Theorem A to the abstract functional differential
equation $\frac{du(t)}{dt}=A(t)u(t)+F(t,u_t)$ with $u_0=\phi\in C$
under appropriate assumptions.

\end{document}